\documentclass[12pt]{article}

\usepackage{amsmath}
\usepackage{amsthm}
\usepackage{amsfonts}
\usepackage{amssymb}
\usepackage{graphicx}

\theoremstyle{plain}
\newtheorem{theorem}{Theorem}[section]
\newtheorem{lemma}{Lemma}[section]
\newtheorem{proposition}{Proposition}[section]

\newtheorem{corollary}{Corollary}[section]

\theoremstyle{definition}
\newtheorem{remark}{Remark}[section]

\title{Chordal Loewner families and univalent Cauchy transforms}
\author{Robert O. Bauer\\ Department of Mathematics\\ University of Illinois at Urbana-Champaign\\ 1409 West Green Street \\ Urbana, IL 61801, USA\\
rbauer@math.uiuc.edu}

\begin{document}

\maketitle

\begin{abstract}
	We study chordal Loewner families in the upper half-plane and 	show that they have a parametric representation. We show 	one, that to every chordal Loewner family there corresponds a 	unique measurable family of probability measures on the real 	line, and two, that to every measurable family of probability 	measures on the real line there corresponds a unique chordal 	Loewner family. In both cases the correspondence is being 	given by solving the chordal Loewner equation. We use this to 	show that any probability measure on the real line with finite 	variance and mean zero has univalent Cauchy transform if and 	only if it belongs to some chordal Loewner family. If the 	probability measure has compact support we give two further 	necessary and sufficient conditions for the univalence of the 	Cauchy transform, the first in terms of the transfinite diameter 	of the complement of the image domain of the reciprocal 	Cauchy transform, and the second in terms of moment 	inequalities corresponding to the Grunsky inequalities.
\end{abstract}

\section{Introduction}

In this paper we discuss chordal Loewner families, the chordal Loewner equation, and probability measures on the real line whose reciprocal Cauchy transform is univalent in the upper half-plane.

Reciprocal Cauchy transforms of probability measures on the real line play an important role in describing the sum of two noncommutative random variables, namely for the free additive convolution developed by Voiculescu \cite{voiculescu:1986}, and the monotonic convolution developed by Muraki \cite{muraki:2000}.

In \cite{schramm:2000}, Schramm introduced a family of random compact sets, growing in a domain of the complex plane. He showed that any random, growing, and compact set that satisfies a certain Markovian-type and conformal invariance property belongs to this family, and that it can be generated by solving Loewner's equation driven by a Brownian motion on the boundary of the domain. This family is now known as stochastic (or Schramm-) Loewner evolution (SLE). Its discovery soon lead to rigorous proofs of various conjectures of conformal field theory about the behavior of certain statistical mechanical systems at criticality, see \cite{werner:2002} and references therein.

In \cite{bauer:2002a}, we noted that a solution of the (chordal) Loewner equation at a fixed time is the reciprocal Cauchy transform of some probability measure on the real line. Since any solution of Loewner's equation takes values in the set of univalent functions this raised the question of what characterizes probability
measures whose reciprocal Cauchy transform is univalent in the upper half-plane. In particular, does any such measure arise by solving a suitable Loewner equation, and if so, what kind of driving functions need to be considered?

To begin to treat this question we found it necessary to study the chordal Loewner equation beyond the cases we found in the literature. These being either to narrow for our purposes, such as the case of compact complement for SLE, \cite{lawler:2001}, or to general, as in \cite{betker:1990}, where, at least to our knowledge, no consistent normalization and thus parametrization of chordal Loewner families with a complete correspondence with driving functions is possible. On the other hand, for the (radial) Loewner equation on the unit disk $\mathbb D$ there exists just such a treatment, given in \cite{rosenblum.rovnyak:1994}. In that case it is convenient to normalize a univalent function $f$ on $\mathbb D$ by $f(0)=0$ and $f'(0)>0$. (Radial) Loewner families, i.e maximal subordination chains of such functions, are then parametrized by the derivative at $z=0$ and one can show that there is a one-to-one correspondence between (radial) Loewner families and so-called Herglotz families, the correspondence being given by solving the (radial) Loewner equation.

 In the chordal case in the upper half-plane we have to deal with compactness questions that do not arise in the (radial) disk case. A suitable class of univalent functions to consider are those  $f$ that map the upper half-plane into the upper half plane and satisfy 
\[
	|f(z)-z|\le \frac{C}{\Im(z)}
\]
for some $C>0$ for all $z$ in the upper half-plane. Such functions are in fact reciprocal Cauchy transforms of probability measures on the real line with finite variance and mean zero. We show that the least constant $C$ in the above inequality serves as a parameter for chordal Loewner families and that chordal Loewner families are in one-to-one correspondence with measurable families of probability measures on the real line, the correspondence being given by solving the chordal Loewner equation. The structure of  our proof of these results is identical to the structure of the proof of the analogous result in the radial case in       
\cite{rosenblum.rovnyak:1994}. However, the basic tools and inequalities used at the various steps in the argument are very different. We give a detailed proof in Sections \ref{S:CLF} and \ref{S:CLE}, taking up the bulk of this paper. We hope that our general treatment may be of use in the context of the stochastic Loewner evolution if the driving function---Brownian motion---is replaced by more general stochastic processes, for example superprocesses. 

As a consequence of our results in Sections \ref{S:CLF} and \ref{S:CLE} we can answer the question, whether
every probability measure on the real line with univalent Cauchy transform belongs to some chordal Loewner family, in the affirmative, at least when the probability measure has finite variance. 

In the case  where the probability measure has compact support we give two further characterizations based on classical results in the theory of univalent functions. The first characterization is in terms of the transfinite diameter of the complement of the image, and is a consequence of a Theorem by Hayman. The second is an application of the Grunsky inequalities. It gives, at least in principle, a characterization of probability measures with univalent Cauchy transform in terms of the moments of the measure.

The paper is structured as follows.
In Section~\ref{S:prelim} we fix notation and collect some results about reciprocal Cauchy transforms of probability measures on the real line. In Section~\ref{S:UCT} we begin by recalling some general results on domains of univalence of reciprocal Cauchy transforms and then obtain three characterizations of univalent Cauchy transforms, Theorems~\ref{T:UCTandCLF} and \ref{T:UCTandCLE}, Corollary~\ref{C:hayman}, and Theorem~\ref{T:grunsky}. In Section~\ref{S:CLF} we introduce and describe chordal Loewner families, culminating in the representation as  parametrized families in Theorem~\ref{T:representation}. Finally, in Section~\ref{S:CLE} we show in Theorem~\ref{T:chordalloewnerequation} that to every chordal Loewner family there corresponds a unique measurable family of probability measures on the real line, where the correspondence is being given by solving the chordal Loewner equation, and in Theorem~\ref{T:measurestoclf} that to every measurable family of probability measures on the real line there corresponds a unique chordal Loewner family, the correspondence again being given by solving the chordal Loewner equation.

The author would like to thank Hari Bercovici for asking the question that inspired this paper.
  
  
\section{Preliminaries}\label{S:prelim}

For the complex plane $\mathbb C$ denote 
$\mathbb H\equiv\{z\in\mathbb C:\Im(z)>0\}$ the upper half-plane, $-\mathbb H\equiv\{z\in\mathbb C:\Im(z)<0\}$ the lower half-plane, and for every positive real number $a$, let 
$\mathbb H_a=\{z\in\mathbb C:\Im(z)>a\}$.
Let $\mu$ be a finite positive Borel measure on $\mathbb R$.
The Cauchy transform $G=G_{\mu}$ of $\mu$ is defined by
\[
	z\in\mathbb H
	\mapsto
	G(z)=\int_{\mathbb R}\frac{\mu(dx)}{z-x}
	\in-\mathbb H.
\]
$G$ is an analytic function with the property
\begin{equation}\label{E:cauchytran}
	\limsup_{y\to\infty}y|G(iy)|<\infty.
\end{equation}
In fact, $\limsup_{y\to\infty}y|G(iy)|=\mu(\mathbb R)$.
Conversely, every analytic function mapping $\mathbb H$ into $-\mathbb H$ that satisfies \eqref{E:cauchytran} is the Cauchy transform of a finite positive Borel measure on $\mathbb R$, \cite[Satz 3, Teil 59, Kapitel VI]{achieser.glasmann:1954}. We can recover $\mu$ from its Cauchy transform  using Stieltjes' inversion formula	
\[
	\mu((a,b))+\mu([a,b])=-\frac{2}{\pi}
	\lim_{\epsilon\searrow0}\int_a^b \Im(G(x+i\epsilon))\ dx.
\]
Since $G(z)\neq0$ for all $z\in\mathbb H$ the reciprocal Cauchy transform $F\equiv1/G$ is an analytic function that maps $\mathbb H$ into $\mathbb H$. Thus $F$ is  a Pick function for which the following representation is known.

\begin{theorem}[Nevanlinna Representation, \cite{bhatia:1997}]
	Every analytic function $F$ such that $\Im (F(z))\ge0$ for 	$z\in\mathbb H$ has a representation
	\[
		F(z)=b+cz+\int_{\mathbb R}\frac{1+tz}{t-z}\ \nu(dt), 
		\quad z\in\mathbb H,
	\]
	where $b=\bar{b}$, $c\ge0$ and $\nu$ is a finite nonnegative 	Borel measure on $\mathbb R$. The triple $(b,c,\nu)$ is unique 	and satisfies $b=\Re(F(i))$,
	 \[
		c=\lim_{0<y\to\infty}\frac{\Im(F(iy))}y,
	\]
	and $\nu(\mathbb R)=\Im(F(i))-c$.
\end{theorem} 
Besides being a Pick function, the reciprocal Cauchy transform $F$ of a probability measure $\mu$  satisfies
\begin{equation}\label{E:decim}
	\inf_{z\in\mathbb H}\frac{\Im(F(z))}{\Im(z)}=1,
\end{equation}
see \cite{maassen:1992}, and the following characterization is known

\begin{theorem}\label{T:decim}\cite{maassen:1992}
	For an analytic function $F:\mathbb H\to\mathbb H$ the 	following are equivalent:
	
	(i) $F$ is the reciprocal Cauchy transform of a probability 	measure $\mu$ on $\mathbb R$. 

	(ii) There exist a real number $b\in\mathbb R$ and a finite 	nonnegative Borel measure $\nu$ on $\mathbb R$ such that 
	\[
		F(z)=b+z+\int_{\mathbb R}\frac{1+tz}{t-z}\ \nu(dt), \quad 
		z\in\mathbb H.
	\]
	
	(iii) $F$ satisfies equation \ref{E:decim}.
\end{theorem}

For probability measures with finite variance and zero mean this result can be specified to

\begin{proposition}\label{P:reciprocal}\cite{maassen:1992}
	For an analytic function $F:\mathbb H\to\mathbb H$ the 	following are equivalent:

	(i) $F$ is the reciprocal Cauchy transform of a probability	 measure on $\mathbb R$ with finite variance and mean zero.

	(ii) There exists a finite positive measure $\rho$ on $\mathbb	R$ such that for all $z\in\mathbb H$,
	\[
		F(z)=z-\int_{\mathbb R}\frac{\rho(dx)}{z-x}.
	\]

	(iii) There exists a positive number $C$ such that  for all 	$z\in\mathbb H$,
	\[
		|F(z)-z|\le\frac{C}{\Im(z)}.
	\]

	Moreover, the variance $\sigma^2$ of $\mu$ in  (i), the total 	weight $\rho(\mathbb R)$  of $\rho$ in (ii), and the  smallest 	possible constant $C$ in (iii) are all equal.
 \end{proposition}

\section{Univalent Cauchy transforms}\label{S:UCT}

We first recall some general results about domains of univalence for Cauchy transforms.   
Denote $\Gamma_{\alpha,\beta}$ the Stolz angle
\[
	\Gamma_{\alpha,\beta}=\{z\in\mathbb H:|z|>\beta
	\text{ and }-\alpha\Im(z)<\Re(z)<\alpha\Im(z)\}.
\]

\begin{proposition}[\cite{bercovici.voiculescu:1993}]
	Let $\mu$ be a probability measure on $\mathbb R$, and let 	$0<\epsilon<\alpha$. There exists a $\beta>0$ such that 

	$(i)$ 	$F=1/G$ is univalent in $\Gamma_{\alpha,\beta}$, and 
	
	$(ii)$ 	$F(\Gamma_{\alpha,\beta})\supset
	\Gamma_{\alpha-\epsilon,\beta(1+\epsilon)}$.
\end{proposition}

For completeness we reproduce the proof of \cite{bercovici.voiculescu:1993}.

\begin{proof}
	Choose $\beta$ so that $|F(z)-z|<\epsilon|z|$ for $z\in 	\Gamma_{\alpha-\epsilon,\beta(1+\epsilon)}$. To prove the 	proposition it is enough to show that each 	$w\in\Gamma_{\alpha-\epsilon,\beta(1+\epsilon)}$ is assumed 	by $F$ exactly once in $\Gamma_{\alpha,\beta}$. Indeed, if 	$\beta'>\beta$ is sufficiently large, then the boundary of 	$\{z\in\Gamma_{\alpha,\beta}:|z|<\beta'\}$ is mapped by $F$ 	into a curve surrounding $w$ exactly once.
\end{proof}

For a probability measure $\mu$ with finite variance there is a stronger result.

\begin{proposition}\label{P:range}
	Let $\mu$ be a probability measure on $\mathbb R$ with finite 	variance $\sigma^2$ and reciprocal Cauchy 	transform $F$. Then the restriction of $F$ to $\mathbb 	H_{\sigma}$ takes every value in $\mathbb H_{2\sigma}$ 	precisely once. 
\end{proposition}

\begin{proof}
	This follows immediately from 
	\cite[Lemma 	2.4]{maassen:1992} 	where the result is 	established under the additional assumption that $\mu$ has 	mean-value 0. Indeed, if $\mu$ has mean value $a$, set 	$\tilde{F}(\cdot)=F(\cdot+a)$. Then $\tilde{F}$ is the reciprocal 	Cauchy transform of $\tilde{\mu}$, where 	$\tilde{\mu}$ is the push-forward of $\mu$ under the 	map $x\mapsto x-a$. $\tilde{\mu}$ has mean value 0 and so 	Lemma 2.4 in \cite{maassen:1992} applies to $\tilde{F}$. This 	in turn implies the result for $F$.
\end{proof}

It follows that there is a right-inverse  $F^{-1}:\mathbb H_{2\sigma}\to\mathbb H_{\sigma}$ and hence that 
$F$ is univalent on $F^{-1}(\mathbb H_{2\sigma})$.   

We now come to the question that was the initial impetus for this paper, namely, 
{\it
when is the Cauchy transform of a probability measure univalent in the entire upper half-plane?}

As a consequence of our general investigation of chordal Loewner families in Section \ref{S:CLF} we have the following

\begin{theorem}\label{T:UCTandCLF}
	Suppose that $\mu$ is a probability measure on  the real line 	with variance $\sigma^2$ and mean zero. The reciprocal 	Cauchy transform $F$ of $\mu$ is univalent in $\mathbb H$ if 	and only if there is a chordal Loewner family 	$\{f(t;\cdot),t\in[0,\infty)\}$ such that $F(z)=f(\sigma^2;z)$, 	$z\in\mathbb H$.
\end{theorem} 

Using the relation between chordal Loewner families and the chordal Loewner equation that we develop in Section \ref{S:CLE} we get

\begin{theorem}\label{T:UCTandCLE}
	Suppose that $\mu$ is a probability measure on  the real line 	with variance $\sigma^2$ and mean zero. The reciprocal 	Cauchy transform $F$ of $\mu$ is univalent in $\mathbb H$ if 	and only if there is a measurable family $\{\mu_t,t\in[0,\infty)\}$
	of probability measures on $\mathbb R$
	such that if we define the family $\{f(t;z),t\in[0,\infty)\}$ as the 	unique solution to the initial value problem
	\[
		\frac{\partial}{\partial t}f(t;z)=-\int_{\mathbb R}
		\frac{\mu_t(dx)}{z-x}\cdot\frac{\partial}{\partial z}f(t;z),
		\quad f(0;z)=z,
	\]
	then $F(z)=f(\sigma^2;z)$, $z\in\mathbb H$.
\end{theorem}

For the definition of chordal Loewner family and the precise meaning of the differential equation see  Sections \ref{S:CLF} and \ref{S:CLE}.

In the case where $\mu$ has compact support we now provide two further characterizations of the univalence of the reciprocal Cauchy transform  $F_{\mu}$. The first characterization is a consequence of a result by Hayman about the transfinite diameter of the ``omitted set'' under meromorphic functions, and the second characterization is in terms of moment conditions based on the Grunsky inequalities. 

We begin by recalling Hayman's result \cite{hayman:1966}.
\begin{theorem}\label{T:Hayman}
	Suppose that $f$ is meromorphic in a domain $D$ whose 	complement $E$ is compact and that $f$ maps $D$ into a 	domain $D'$ whose complement is $E'$. Further suppose that 	$f'(\infty)=1$ which means that 
	\[
		f(z)=z+a_0+\frac{a_1}{z}+\cdots\quad \text{for large }z.
	\]
	Then $d(E')\le d(E)$, where $d(E)$ and $d(E')$ denote the 	transfinite diameter of $E$ and $E'$ respectively. Equality 	holds if $f$ is univalent and maps $D$ onto $D'$.
\end{theorem}

Let $\mu$ be a compactly supported probability measure on the real line with Cauchy transform $G=G_{\mu}$ and reciprocal Cauchy transform $F=F_{\mu}$. Denote $[A_{\mu}, B_{\mu}]$ the convex closure of the support of $\mu$. Using the Schwarz reflection principle it is easy to see that both $G$ and $F$ extend as analytic functions to $\mathbb C\backslash[A_{\mu}, B_{\mu}]$, \cite{bauer:2003a}. Denote these extensions also by $G$ and $F$. Then we have the following
 
\begin{corollary}\label{C:hayman}
	With the notation from above, $\mu$ has Cauchy transform 	univalent in the upper half-plane if and only if the transfinite 	diameter of the complement of $F(\mathbb C 	\backslash[A_{\mu}, B_{\mu}])$ equals $B_{\mu}-A_{\mu}$.
\end{corollary}

Next we consider the characterization by moments. To simplify notation we will assume that the support of $\mu$ is contained in the interval $[-2,2]$. Let  
\[
	a_n=\int_{\mathbb R}x^n\ \mu(dx),\quad n=0,1,2,\dots,
\]
and note that
\[
	G(z)=\sum_{n=0}^{\infty}a_n z^{-(1+n)},\quad |z|>2.
\]
$G$ extends as an analytic function to $\mathbb C\backslash[-2,2]$. Define 
\[
	\psi:\{z\in \mathbb C:|z|>1\}\to\mathbb C\backslash[-2,2]
\]
by $\psi(z)=z+1/z$. Then $\psi$ is univalent and onto. From the expansion
\[
	\left(z+\frac{1}{z}\right)^{-1}=\frac{1}{z}\cdot\frac{1}{1+z^{-2}}
	=\sum_{k=0}^{\infty}(-1)^k z^{-(2k+1)}
\]
valid for $|z|>1$, we get
\[
	\left(z+\frac{1}{z}\right)^{-(n+1)}=\sum_{k=0}^{\infty}
	(-1)^k\binom{n+k}{n}z^{-(n+2k+1)},\quad|z|>1,
\]
since $\binom{n+k}{n}$ is the number of distinct, positive, and odd integer solutions to $x_1+\cdots+x_{n+1}=n+2k+1$. Thus, by rearranging,
\begin{align}
	G(\psi(z))&=\sum_{n=0}^{\infty}a_n \sum_{k=0}^{\infty}
	(-1)^k\binom{n+k}{n}z^{-(n+2k+1)}\notag\\
	&=\sum_{n=0}^{\infty}\left[\sum_{k=0}^{[n/2]}a_{n-2k}(-1)^k
	\binom{n-k}{n-2k}\right]\ z^{-(n+1)},\notag
\end{align}
and the latter expansion holds for $|z|>1$. For $n=0,1,2,\dots$, set
\[
	\alpha_n=\sum_{k=0}^{[n/2]}a_{n-2k}(-1)^k
	\binom{n-k}{n-2k}.
\]
If $F=F_{\mu}=1/G$, then
\[
	F(\psi(z))=\frac{1}{G(\psi(z))}=\sum_{n=0}^{\infty}\beta_n 	z^{1-n},
\]
where
\begin{align}
	\alpha_0\beta_0&=1,\notag\\
	\alpha_0\beta_1+\alpha_1\beta_0&=0,\notag\\
	\alpha_0\beta_2+\alpha_1\beta_1+\alpha_2\beta_0&=0,
	\notag\\ 
	\dots&\notag
\end{align}
Since $a_0=\mu(\mathbb R)=1$, we have $\alpha_0=1$ and $\beta_0=1$. Solving the above system inductively and substituting back the $a_n$s for the $\alpha_n$s we find for instance $\beta_1=-a_1$, $\beta_2=1+a_1^2-a_2$, and 
$\beta_3=-a_1^3+2a_1 a_2 -a_3$. For the function $F\circ\psi$ we can now consider the Grunsky inequalities. 

We briefly recall the definition for the Grunsky coefficients. All results we use regarding these coefficients can be found in \cite{duren:1983}. For an analytic function 
$g$  with an expansion
\[
	g(z)=z+b_0+b_1 z^{-1}+b_2 z^{-2}+\cdots
\]
valid for $|z|>1$, consider
\[
	\frac{\zeta g'(\zeta)}{g(\zeta)-w}=\sum_{n=0}^{\infty} 	F_n(w)\zeta^{-n},
\]
where the expansion is valid for all $\zeta$ in some neighborhood of $\infty$ with $F_n(w)=w^n+\sum_{k=1}^n a_{nk}w^{n-k}$, the $n$-th Faber polynomial of $g$. Then
\[
	F_n(g(z))=z^n+\sum_{k=1}^{\infty}\beta_{nk}z^{-k},\quad n=1,
	2,\dots
\]
The coefficients $\beta_{nk}$ are known as the Grunsky coefficients of $g$. Set 
\[
	c_{nk}=\sqrt{\frac{k}{n}}\beta_{nk},\quad 
	(n,k)\in(\mathbb Z^+)^2.
\]
Then the (weak) Grunsky inequalities hold if for each $N\in\mathbb Z^+$, $(\lambda_1,\dots,\lambda_N)\in\mathbb C^N$,
\[
	\left|\sum_{n=1}^N\sum_{k=1}^N 	c_{nk}\lambda_n\lambda_k\right|  	\le\sum_{n=1}^N|\lambda_n|^2.
\]
These inequalities are a necessary and sufficient condition for $g$ to be univalent on $\{z\in\mathbb C:|z|>1\}$. We now apply this fact  to $F\circ\psi$.

\begin{theorem}\label{T:grunsky}
	Suppose that $\mu$ is a probability measure on $\mathbb R$ 	such that its support is contained in $[-2,2]$. Then the 	reciprocal Cauchy transform $F$ of $\mu$ is univalent in 	$\mathbb H$ if and only if for each $N\in\mathbb Z^+$ the 	real symmetric matrix $[c_{nk}]_{n,k=1}^N$ has all its 	eigenvalues in $[-1,1]$.
\end{theorem}

\begin{proof}
	Apply the Grunsky inequalities to $F\circ\psi$. Since all 	coefficients of $F\circ\psi$ are real, its Grunsky coefficients 	are also real. Now the Grunsky inequalities reduce to  bounds 	on the eigenvalues of the matrices $[c_{nk}]_{n,k=1}^N$.
\end{proof}

These conditions become quickly intractable. For $N=1$ we get
$|c_{11}|=|1+a_1^2-a_2|\le1$, while for $N=2$ the matrix 
$[c_{nk}]_{n,k=1}^2$ reads
\begin{equation}
	\left[
	\begin{matrix}
		1+a_1^2-a_2 &-\sqrt{2}(a_1^3-2a_1 a_2+a_3)\\
		-\sqrt{2}(a_1^3-2a_1 a_2+a_3) & 
		1+3a_1^4-8a_1^2a_2+3a_2^2+4a_1a_3-2a_4\notag
	\end{matrix}
	\right].
\end{equation}
If $\mu$ is even, i.e. all odd moments vanish, then this gives the two conditions, $a_2\le2$ and  $|1+3a_2^2-2a_4|\le 1$. Thus the Grunsky inequalities may be useful to quickly rule out that a certain distribution has a reciprocal Cauchy transform univalent in the upper half-plane just by looking at a few of its moments. 	 

\section{Chordal Loewner families}\label{S:CLF}

Denote $\mathfrak R$ the class of analytic functions $f:\mathbb H\to\mathbb H$ which are univalent and satisfy
\begin{equation}\label{E:keybound}
	|f(z)-z|\le\frac{C}{\Im(z)},\quad z\in\mathbb H,
\end{equation}
for some constant $C\in[0,\infty)$. Denote $a$ the least such constant. By Proposition \ref{P:reciprocal}, part (ii),
\[
	iy[iy-f(iy)]=\int_{\mathbb R}\frac{iy}{iy-x}\ \rho(dx),
\]
where $\rho$ is a nonnegative Borel measure on $\mathbb R$ with total mass $a$. Thus, by bounded convergence,
\[
	\lim_{y\to\infty}iy[iy-f(iy)]=\rho(\mathbb R)=a.
\]

\begin{remark}\label{R:compact}
	If $K\subset\overline{\mathbb H}$ is compact and such that 
	$\mathbb H\backslash K$ is connected and simply 	connected, then there exists a unique $f\in\mathfrak R$ such 	that 
	$f(\mathbb H)=\mathbb H\backslash K$, see 	\cite{lawler:2001}.
\end{remark}

For an example, let $t\in[0,\infty)$ and define the function $f(t;\cdot)$ by
\[
	f(t;z)=\sqrt{z^2-2t}
	=z-\frac{t}{z}+O\left(\frac{1}{|z|^2}\right),\quad z\to\infty.
\]
Then $f(t;\cdot)$ belongs to $\mathfrak R$ and its range $G_t$ is the upper half-plane with a slit along the imaginary axis from zero to $\sqrt{2t}$. We note that $f(t;\cdot)$ is the reciprocal Cauchy transform of the arcsine law with density $1/(\pi\sqrt{2t-x^2})$ supported in $[-\sqrt{2t}, \sqrt{2t}]$. The functions $f(t;\cdot)$, $t\in[0,\infty)$, form a totally ordered ``chain'' relative to the partial ordering induced by inclusion of the image domains. In fact, $\{G_t\}_{t\in[0,\infty)}$ is a maximal totally ordered family of simply connected regions in $\mathbb H$.   

Let $f, g\in\mathfrak R$. We say $f$ is subordinate to $g$ and write $f\prec g$ if $f=g\circ h$ for some $h\in\mathfrak R$. 

\begin{lemma}\label{L:subordinated}
	If $f,g\in\mathfrak R$, then $f\prec g$ if and only if 
	$f(\mathbb H) \subseteq g(\mathbb H)$. In this case
	\begin{equation}\label{E:limincrease}
		\lim_{y\to\infty}iy\left[iy-f(iy)\right]
		\le\lim_{y\to\infty}iy\left[iy-g(iy)\right]
	\end{equation}
	with equality if and only if $f\equiv g$.
\end{lemma}

\begin{proof}
	If $f=g\circ h$ with $h\in\mathfrak R$, then 	
	$f(\mathbb H)=g(h(\mathbb H))\subseteq g(\mathbb H)$, 	so the condition is necessary. Conversely, if 
	$f(\mathbb H) \subseteq g(\mathbb H)$, then 
	$h\equiv g^{-1}\circ f:\mathbb H\to\mathbb H$ is univalent. 	Denote $a,b\in(0,\infty)$ the least constants such that
	\begin{equation}\label{E:fbound}
		|f(z)-z|\le\frac{b}{\Im(z)},\quad z\in\mathbb H,
	\end{equation}
	\begin{equation}\label{E:gbound}
		|g(z)-z|\le\frac{a}{\Im(z)}\quad z\in\mathbb H.
	\end{equation}
	We then also have 
	\begin{equation}\label{E:flim}
		f(iy)=i\left(y+\frac{b}{y}\right)
		+o\left(\frac{1}{y}\right),\quad y\to\infty,
	\end{equation}
	\begin{equation}\label{E:glim}
		g(iy)=i\left(y+\frac{a}{y}\right)
		+o\left(\frac{1}{y}\right),\quad y\to\infty.
	\end{equation}
	By \eqref{E:flim} and \eqref{E:gbound}, 
	$f(iy)\in g(\mathbb H_1)$ for all $y$ large enough, and if 	$f(iy)=g(z)$ for some $z$ with $\Im(z)>1$, then
	\[
		|g^{-1}(f(iy))-f(iy)|=|z-g(z)|<a.
	\]
	Since also $|f(iy)-iy|<b$ for $y>1$ we get
	\[
		|g^{-1}(f(iy))-iy|\le|g^{-1}(f(iy))-f(iy)|+|f(iy)-iy|<a+b
	\]
	for all $y$ large enough. Hence
	\[
		\lim_{y\to\infty}\frac{y}{|g^{-1}(f(iy))|}=1
	\]
	and $g^{-1}\circ f$ is the reciprocal Cauchy transform of a 	probability measure on the real line. Let $z=g^{-1}(iy)$ and set 	$\tilde{y}=y-a/y$. For $y$ large enough both $\Im(z)\ge 1$ and 	$\tilde{y}\ge 1$. Then $g(i\tilde{y})=iy+o(1/|y|)$ and so
	\[
		|z-i\tilde{y}|=|g^{-1}(iy)-g^{-1}(g(i\tilde{y}))| 		=o\left(\frac{1}{|y|}\right)
	\]
	since $|(g^{-1})'(z)|$ is bounded for $\Im(z)\ge 1$. Similarly 	$|g'(z)|$ is bounded for $\Im(z)\ge1$ and so
	\[
		|g(z)-g(i\tilde{y})|\le C|z-i\tilde{y}|=o\left(\frac{1}{|y|}\right).
	\]
	In particular $z=O(|y|)$ and
	\begin{align}
		iy[iy-g^{-1}(iy)]&=g(z)[g(z)-z]\notag\\
		&=z[g(z)-z]+o(1)=i\tilde{y}[g(i\tilde{y})-i\tilde{y}]+o(1).\notag
	\end{align}
	So
	\[
		\lim_{y\to\infty}iy[iy-g^{-1}(iy)]
		=\lim_{\tilde{y}\to\infty}i\tilde{y}[g(i\tilde{y})-i\tilde{y}]=-a.
	\]
	Since
	\begin{align}
		iy[iy-g^{-1}(f(iy))]&=iy[iy-f(iy)]+iy[f(iy)-g^{-1}(f(iy))]\notag\\
		&=iy[iy-f(iy)]\notag\\
		&\ +i\left(y+\frac{b}{y}\right)\left[i\left(y+\frac{b}{y}\right)
		-g^{-1}\left(i\left(y+\frac{b}{y}\right)\right)\right]+o(1),\notag
 	\end{align}
	we now get
	\begin{equation}\label{E:ginversef}
		\lim_{y\to\infty}iy[iy-g^{-1}(f(iy))]=b-a.
	\end{equation}
	If $F$ is the reciprocal Cauchy transform of a probability 	measure $\mu$ on $\mathbb R$ we introduce the function
	\[
		y\in(0,\infty)\mapsto C_F(y)\equiv 			y\left(\frac{1}{F(iy)}-\frac{1}{iy}\right)\in\mathbb C.
	\]
	One can show, \cite{maassen:1992}, that 
	\[
		\lim_{y\to\infty}y\Im(C_F(y))=\int_{\mathbb R}x^2\ \mu(dx),
	\]
	and, if $\int_{\mathbb R}x^2\ \mu(dx)<\infty$, then 
	\[
		\lim_{y\to\infty}\Re(C_F(y))=\int_{\mathbb R}x\ \mu(dx).
	\]
	On the other hand,
	\[
		C_F(y)=-\frac{iy}{F(iy)}[iy-F(iy)] 
	\]
	and so for $F=g^{-1}\circ f$
	\[
		\lim_{y\to\infty}y C_F(y)
		=-\lim_{y\to\infty}\frac{y}{F(iy)}\lim_{y\to\infty}iy[iy-F(iy)]
		=i(b-a).
	\]
	Hence $g^{-1}\circ f$ is the reciprocal Cauchy transform of a 	probability measure on the real line with mean zero and 	variance $b-a$. It follows in particular that $f\prec g$ and also 	that $b\ge a$, i.e. \eqref{E:limincrease}. Finally, if $b=a$, then 	$\mu$ has mean and variance both equal to zero, that is 	$\mu$ is a unit point mass at zero. In that case $\int_{\mathbb 	R}\mu(dx)/(z-x)=1/z$ and so $g^{-1}(f(z))=z$. 
\end{proof}

Note that the above result implies in particular that the functions $f$ in $\cal R$ are determined by their image domains.
Subordination is a partial ordering on $\mathfrak R$. By a chain in $\mathfrak R$ we mean a nonempty totally ordered subset $\mathfrak C$ of $\mathfrak R$, and by a {\it chordal Loewner family} $\mathfrak L$ we mean a maximal chain, i.e. whenever ${\mathfrak L}\subseteq\mathfrak C$ where $\mathfrak C$ is a chain in $\mathfrak R$, then $\mathfrak L=\mathfrak C$.

\begin{theorem} Every $f\in\mathfrak R$ belongs to some chordal 	Loewner family. More generally, every chain $\mathfrak C$ in 	$\mathfrak R$ is contained in a chordal Loewner family.
\end{theorem}

\begin{proof}
	See proof of Theorem B in section 7.10 in \cite{rosenblum.rovnyak:1994}.
\end{proof}

The following theorem plays in the chordal case the role the Carath\'eodory convergence theorem plays in the radial (disk) case.

\begin{theorem}\label{T:caratheodory}
	Let $\{f_n\}_{n=1}^{\infty}$ be a sequence in $\mathfrak R$ 	and for every $n\in\mathbb Z^+$ let 
	$a_n=\lim_{y\to\infty}iy[iy-f_n(iy)]$. If $a=\lim_{n\to\infty}a_n$ 	exists and is finite, then there exists a function $f\in\mathfrak 	R$ 	and a subsequence $n_1, n_2,\dots$, such 	that for every $m\in\mathbb Z^+$
	\begin{equation}\label{E:uniformlim}
		\sup_{z\in\mathbb H_{1/m}}|f(z)-f_{n_k}(z)|\to0,\quad
		\text{as } k\to\infty.
	\end{equation}
	Furthermore, if $A^o$ denotes the interior of a set $A$, then
	\begin{equation}\label{E:imagedomain}
		f(\mathbb H)=\bigcup_{m=1}^{\infty}\bigcup_{l=1}^{\infty}
		\left(\bigcap_{k=l}^{\infty}f_{n_k}(
		\mathbb H_{1/m})\right)^o,
	\end{equation}
	and
	\begin{equation}\label{E:limita}
		a=\lim_{y\to\infty} iy[iy-f(iy)].
	\end{equation}
	If, in addition, the family $\{f_n\}_{n=1}^{\infty}$ is totally 	ordered, then $f$ is unique and \eqref{E:uniformlim} and 	\eqref{E:imagedomain} hold without going to subsequences.
\end{theorem}  

\begin{proof}
	Since
	\[
		|f_n(z)|=|(f_n(z)-z)+z|\le\frac{a_n}{\Im(z)}+|z|,
	\]
	the family $\{f_n\}$ is locally bounded and there exists a 	subsequence $n_1, n_2,\dots$ and an analytic function 	$f:\mathbb H\to\mathbb H$ such that $f_{n_k}\to f$ uniformly 	on compact subsets of $\mathbb H$ as $k\to\infty$. 	Furthermore, since $f_n$ is univalent for each $n$, $f$ is either 	univalent or constant. For each $z\in\mathbb H$,
	\begin{align}
		|f(z)-z|&=|(f(z)-f_{n_k}(z))-(f_{n_k}(z)-z)|\notag\\
		&\le|f(z)-f_{n_k}(z)|+\frac{a_{n_k}}{\Im(z)}.\notag
	\end{align}
	It follows that $f$ is nonconstant, that $|f(z)-z|\le a/\Im(z)$, and 	that $a$ is the least constant such that this inequality holds. 	This proves \eqref{E:limita} along a subsequence.
	
	We have
	\[
		f_n(z)=z-\int_{\mathbb R}\frac{\rho_n(dx)}{z-x},\quad  
		f(z)=z-\int_{\mathbb R}\frac{\rho(dx)}{z-x},
	\]
	where $\rho_n$ and $\rho$ are nonnegative Borel measures 	with total mass $a_n$ and $a$, respectively. Suppose that 	$a>0$. Then 
	\[
		\frac{1}{a_{n_k}}\int_{\mathbb R}\frac{\rho_{n_k}(dx)}{z-x}
		\to\frac{1}{a}\int_{\mathbb R}\frac{\rho(dx)}{z-x},
	\]
	uniformly on compact subsets of $\mathbb H$. Since 	$\rho_n/a_n$ and $\rho/a$ are probability measures, it follows 	by \cite[Theorem 2.5]{maassen:1992} that $\rho_n/a_n$ 	converges weakly to $\rho/a$. This implies that 
	\[
		\int_{\mathbb R}\frac{\rho_{n_k}(dx)}{z-x}
		\to\int_{\mathbb R}\frac{\rho(dx)}{z-x}
	\]
	uniformly on $\mathbb H_{1/m}$, for every 
	$m\in\mathbb Z^+$, see \cite{bauer:2003a}. The case $a=0$ 	is easily treated directly. This proves \eqref{E:uniformlim} 	along a subsequence.

	Next, if 
	\[
		w\in\bigcup_{m=1}^{\infty}\bigcup_{l=1}^{\infty}\left(  		\bigcap_{k=l}^{\infty}f_{n_k}(\mathbb H_{1/m})\right)^o,
	\]
	then there exist $m,k\in\mathbb Z^+$ and $\epsilon>0$ such 	that
	\[
		\{w':|w-w'|<\epsilon\}\subset  \left(  		\bigcap_{k=l}^{\infty}f_{n_k}(\mathbb H_{1/m})\right)^o.
	\]
	Consider $A_{n_k}\equiv f_{n_k}^{-1}(\{w':|w-w'|<\epsilon\})
	\subset\mathbb H_{1/m}$. Then, by \eqref{E:uniformlim}, 
	$w\in f(A_{n_k})$ for all $k$ large enough. Conversely,	suppose $w\in f(\mathbb H)$. Then there exist 
	$m\in\mathbb Z^+$ and $\epsilon>0$ such that 
	$\{w':|w-w'|<\epsilon\}\subset f(\mathbb H_{1/m})$. Consider
	$A\equiv f^{-1}(\{w':|w-w'|<\epsilon\})$. Then, by 	\eqref{E:uniformlim}, 
	$f_{n_k}(\mathbb H_{1/m})\supset\{w':|w-w'|<\epsilon/2\}$ for
	all $k$ large enough and this proves \eqref{E:imagedomain} 
	along a subsequence.

	Finally, suppose that $\{f_n\}_{n=1}^{\infty}$ is totally ordered. 
	Given $m,n\in\mathbb Z^+$ we may assume without loss of 	generality that $f_n\prec f_m$. Then $f_n=f_m\circ g$ for 	some $g\in\mathfrak R$. By Lemma \ref{L:subordinated}
	\[
		|g(z)-z|\le\frac{a_n-a_m}{\Im(z)}.
	\]
	Thus $|f_n(z)-f_m(z)|=|f_m(g(z))-f_m(z)|\to 0$, uniformly on 	compact subsets of $\mathbb H$, as $m,n\to\infty$, and this 	proves the theorem.	  
\end{proof}

\begin{lemma}\label{L:subordinatelim}
	Let $f, f_1, f_2,\dots$ and $g, g_1, g_2,\dots$ belong to 
	$\mathfrak R$. For each $n\in\mathbb Z^+$ let 	$a_n=\lim_{y\to\infty}iy[iy-f_n(iy)]$, 
	$b_n=\lim_{y\to\infty}iy[iy-g_n(iy)]$. Assume that $\sup_n 	a_n<\infty$, $\sup_n b_n<\infty$, and $f_n\to f$, $g_n\to g$, 	uniformly on compact subsets of $\mathbb H$, as $n\to\infty$. 	If $f_n\prec g_n$ for each $n\in\mathbb Z^+$, then $f\prec g$.
\end{lemma} 

\begin{proof}
	For each $n\in\mathbb Z^+$, define the function $h_n$ by 	$f_n=g_n\circ h_n$. Then $\lim_{y\to\infty}iy[iy-h(iy)]=a_n-b_n$.
	By Theorem \ref{T:caratheodory}, there exists 
	$h\in\mathfrak R$ and a subsequence $n_1,n_2,\dots$ such 	that $h_{n_k}\to h$ uniformly on compact subsets of 
	$\mathbb H$, as $k\to\infty$.
	Since
	\begin{align}
		|g(h(z))-f(z)|&\le |g(h(z))-g(h_{n_k}(z))|+|g(h_{n_k}(z))-
		g_{n_k}(h_{n_k}(z))|\notag\\
		&\quad+|f_{n_k}(z)-f(z)|,\notag
	\end{align}
	and $g_{n_k}\to g$, $f_{n_k}\to f$ uniformly, the right hand 	side tends to zero as $k\to\infty$.
\end{proof}

\begin{lemma}\label{L:jordanarc}
	Let $\gamma:[0,1)\to\overline{\mathbb H}$ be a Jordan arc 	such that
	 \[
		\gamma(0)\in\partial\mathbb H\quad\text{and}\quad 
		\lim_{t\nearrow1}\Im(\gamma(t))=\infty.
	\]
	 For each $t\in[0,1)$ 	let $f(t;\cdot)\in\mathfrak R$ be the 	unique function whose range is the complement of 	$\gamma([0,t])$ in $\mathbb H$, and set 	$a(t)=\lim_{y\to\infty}iy[iy-f(t;iy)]$. Then 
	$t\in[0,1)\mapsto a(t)\in[0,\infty)$ is nondecreasing, continuous, 	and onto.
\end{lemma}   

\begin{proof}
	Let $t, t_1, t_2,\dots$ be points in $[0,1)$ such that $t_n\to t$, 	as $n\to\infty$. By Remark \ref{R:compact}, $\sup_n 	a(t_n)<\infty$, and so there is a convergent subsequence. 	Applying the first part of Theorem \ref{T:caratheodory}, it 	follows that there is a subsequence $n_1,n_2,\dots$ such that 
	$f(t_{n_k})\to f\in\mathfrak R$, and $a(t_{n_k})\to a$. It is 	straightforward to check that    
	\[
		\bigcup_{m=1}^{\infty}\bigcup_{l=1}^{\infty}
		\left(\bigcap_{k=l}^{\infty}f(t_{n_k};\mathbb H_{1/m})
		\right )^o=\mathbb H\backslash\gamma([0,t]),
	\]
	and so $f=f(t)$. Now we apply the second part of Theorem 
	\ref{T:caratheodory} and it follows that $t\mapsto a(t)$ is 	continuous. It remains to check that $a(t)\to\infty$ as 	$t\nearrow1$. By Proposition \ref{P:range}, $f(t;\mathbb H)$ 	contains $\mathbb H_{2\sqrt{a(t)}}$ and the lemma follows.
\end{proof}

\begin{lemma}\label{L:complete}
	Let $f\prec g$ where $f,g\in\mathfrak R$, 
	$a=\lim_{y\to\infty}iy[iy-f(iy)]$, $b=\lim_{y\to\infty}iy[iy-g(iy)]$, 	and let $c$ be a positive number.

	(i) If $c<b$, there is an $h\in\mathfrak R$ such that 	\[
		\lim_{y\to\infty}iy[iy-h(iy)]=c\quad\text{and}\quad g\prec h.
	\]
	(ii) If $b< c< a$, there is an $h\in\mathfrak R$ such that 	\[
		\lim_{y\to\infty}iy[iy-h(iy)]=c\quad\text{and}\quad f\prec 		h\prec g.
	\]
	(iii) If $a<c$, there is an $h\in\mathfrak R$ such that 	\[
		\lim_{y\to\infty}iy[iy-h(iy)]=c\quad\text{and}\quad h\prec f.
	\]
\end{lemma} 

\begin{proof}
	The result can be reduced to the case where the complements  	of the ranges of $f$ and $g$ in $\mathbb H$ are compact, 	bounded by Jordan arcs, and separated by at least 	$\epsilon>0$ in $\mathbb H$. To reduce to the case of 	compact complement, let $f_n$ be the element of 
	$\mathfrak R$ such that $f_n(\mathbb H)=f(\mathbb H)\cup
	\{z\in\mathbb H:|\Re(z)|>n\}$. By Theorem 	\ref{T:caratheodory}, $f_n\to f$, uniformly on compact subsets 	of $\mathbb H$. Define $g_n$ similarly. By construction, 	$f_n\prec g_n$, and it follows from Lemma 	\ref{L:subordinatelim}
	that it is enough to proof the result for $f_n$ and $g_n$. 	Approximating the compact complement $K$ of the range of 	$f$ in $\mathbb H$  by lemniscates we may assume that $K$ 	is bounded by Jordan arcs, and, after shifting the range of $f$  	by 	a small amount along the imaginary axis, we may assume 	that the complements of the ranges of $f$ and $g$ in 	$\mathbb H$ 	are separated by at least $\epsilon$. 

	Now the proof proceeds as in 
	\cite[Lemma 7.11D]{rosenblum.rovnyak:1994}. The 	Jordan arc $\gamma$ used to produce $h$ first traces out the 	boundary 	of the range of $g$ in $\mathbb H$, say from left to 	right. If this part of the boundary consists of more than one 	component, then the Jordan arc connects the components by 	moving along the real axis between components. After 	$\gamma$ has traced the boundary of $g(\mathbb H)$ in 	$\mathbb H$, it then moves out to trace the boundary of 	$f(\mathbb H)$ in $\mathbb H$ from right to left. Since the 	two boundaries are separated, $\gamma$ continues to be a 	Jordan arc. Finally, after $\gamma$ has traced both 	boundaries it continues to $\infty$ so that its imaginary part 	also goes to $\infty$.            
\end{proof}

\begin{theorem}\label{T:representation}
	If $\mathfrak L$ is any chordal Loewner family, then 
	\[
		f\in{\mathfrak L}\mapsto
		\lim_{y\to\infty}iy[iy-f(iy)]\in[0,\infty) 
	\]
	is one-to-one and onto. Thus the family $\mathfrak L$ has a 	parametric representation 
	$\mathfrak L=\{f(t;\cdot)\}_{t\in[0,\infty)}$, where
	each $f(t;\cdot)$ satisfies
	\[
		f(t;iy)=i\left(y+\frac{t}{y}\right)+o\left(\frac{1}{|y|}\right),\quad
		y\to\infty,
	\]
	$f(b;\cdot)$ is subordinate to $f(a;\cdot)$, whenever $0\le 	a\le b<\infty$, and $f(0;z)=z$.
\end{theorem}

\begin{proof}
	The proof is analogous to the proof of 	\cite[Theorem 7.12]{rosenblum.rovnyak:1994} and is omitted.
\end{proof}

 
\section{Chordal Loewner equation}\label{S:CLE}

Let $\mathfrak L$ be any chordal Loewner family with parametric representation $\{f(t;\cdot),t\in[0,\infty)\}$. If $0\le a\le b<\infty$, then $f(b;\cdot)\prec f(a;\cdot)$ and therefore 
\[
	f(b;z)=f(a;B(a,b;z))
\]
for some function $B(a,b;\cdot)\in\mathfrak R$. Then $B(a,a;z)=z$,
\[
	\lim_{y\to\infty}iy[iy-B(a,b;iy)]=b-a,
\]
and 
\begin{equation}\label{E:Bflow}
	B(a,c;z)=B(a,b;B(b,c;z)),
\end{equation}
whenever $0\le a\le b\le c<\infty$. We say $\mathfrak L$ is a chordal Loewner family with associated semigroup $\{B(a,b;\cdot),0\le a\le b<\infty\}$. Since $f(t;z)=B(0,t;z)$, $t\in[0,\infty), z\in\mathbb H$, the semigroup determines $\mathfrak L$.

\begin{theorem}\label{T:absolute}
	Let $\mathfrak L$ be a chordal Loewner family with associated 	semigroup $\{B(a,b;\cdot),0\le a\le b<\infty\}$. Then for all 	$z\in\mathbb H$,
	\begin{equation}\label{E:absolutefirst}
		|B(a,c;z)-B(b,c;z)|\le\frac{b-a}{\Im(z)},	\end{equation}
	\begin{equation}\label{E:absolutesecond}
		|B(a,b;z)-B(a,c;z)|\le\left(1+\frac{b-a}{\Im(z)^2}\right)
		\frac{c-b}{\Im(z)},
	\end{equation}
	whenever $0\le a\le b\le c<\infty$. Thus, for each 
	$z\in\mathbb H$,
	
	(i) the function $t\in[0,\infty)\mapsto f(t;z)\in\mathbb H$ is 
	absolutely continuous,

	(ii) if $b>0$, $a\in[0,b]\mapsto B(a,b;z)\in\mathbb H$ is 	absolutely continuous,

	(iii) if $a\ge0$, $b\in[a,\infty)\mapsto B(a,b;z)\in\mathbb H$ is 
	absolutely continuous.
\end{theorem}

\begin{proof}
	We have
	\begin{align}
		|B(a,c;z)-B(b,c;z)|&=|B(a,b;B(b,c;z))-B(b,c;z)|\notag\\
		&\le\frac{b-a}{\Im(B(b,c;z))}\le\frac{b-a}{\Im(z)},\notag
	\end{align}
	which proves \eqref{E:absolutefirst}. For 	\eqref{E:absolutesecond}, note that 
	\[
		B(a,b;z)=z-\int_{\mathbb R}\frac{\rho_{a,b}(dx)}{z-x},
	\]
	for some nonnegative Borel measure $\rho_{a,b}$ with 	$\rho_{a,b}(\mathbb R)=b-a$. So
	\[
		|B'(a,b;z)-1|=\left|\int_{\mathbb R}
		\frac{\rho_{a,b}(dx)}{(z-x)^2}\right|\le\frac{b-a}{\Im(z)^2}.
	\]
	Since also $|B(b,c;z)-z|\le(c-b)/\Im(z)$, we get
	\begin{align}
		|B(a,b;z)-B(a,c;z)|&=|B(a,b;z)-B(a,b;B(b,c;z))|\notag\\
		&\le\left(1+\frac{b-a}{\Im(z)^2}\right)
		\frac{c-b}{\Im(z)}.\notag
	\end{align}
\end{proof}

\begin{theorem}\label{T:uniformderivative}
	Assume the same situation as in Theorem \ref{T:absolute}.

	(i) There is a subset $N$ of $[0,\infty)$ of Lebesgue measure 	zero such that if $t\in[0,\infty)\backslash N$, then
	\[
		\frac{\partial}{\partial t}f(t;z)
		=\lim_{h\to0}\frac{f(t+h;z)-f(t;z)}{h}
	\]
	exists, uniformly on compact subsets of $\mathbb H$.

	(ii) For each $b>0$, there is a subset $N$ of $[0,b]$ of 	Lebesgue measure zero such that if $a\in[0,b]\backslash N$, 	then 
	\[
		\frac{\partial}{\partial a}B(a,b;z)
		=\lim_{h\to0}\frac{B(a+h,b;z)-B(a,b;z)}{h}
	\]
	exists, uniformly on compact subsets of $\mathbb H$.

	(iii) For each $a\ge0$, there is a subset $N$ of $[a,\infty)$ of 	Lebesgue measure zero such that if 
	$b\in[a,\infty)\backslash N$, then 
	\[
		\frac{\partial}{\partial b}B(a,b;z)
		=\lim_{h\to0}\frac{B(a,b+h;z)-B(a,b;z)}{h}
	\]
	exists, uniformly on compact subsets of $\mathbb H$.
\end{theorem}

\begin{proof}
	We will only check (i). The other parts can then be handled in a 	similar way. By Theorem \ref{T:absolute}, for fixed $z$, 	$(\partial/\partial t)f(t;z)$ exists $a.e.$ on $(0,\infty)$. The 	exceptional null set depends on $z$, but we may choose a 
	single nullset $N\subset[0,\infty)$ such that the 	derivative exists for all $t\in[0,\infty)\backslash N$ and 
	$z=1/2, 2/3, 3/4,\dots$. Fix $t\in[0,\infty)\backslash N$, and 	consider the difference quotients
	\begin{equation}\label{E:differencequotients}
		\left\{\frac{f(t+h;z)-f(t;z)}{h},0<|h|<t/2\right\}
	\end{equation}
	as analytic functions on $\mathbb H$. Note that 	$f(t+h;z)-f(t;z)=B(0,t+h;z)-B(0,t;z)$, and so, by 	\eqref{E:absolutesecond},
	\begin{equation}
		|f(t+h;z)-f(t;z)|\le
		\begin{cases}
			\left(1+\frac{t}{\Im(z)^2}\right)\frac{h}{\Im(z)}, 
			&\text{if $t/2>h>0$,}\\
			\left(1+\frac{t+h}{\Im(z)^2}\right)\frac{-h}{\Im(z)},
	  		&\text{if $-t/2<h<0$.}\\
		\end{cases}
	\end{equation}
	Thus, for all $0<|h|<t/2$,
	\[
		\left|\frac{f(t+h;z)-f(t;z)}{h}\right|
		\le\left(1+\frac{t}{\Im(z)^2}\right)\frac{1}{\Im(z)}
	\]
	and the family \eqref{E:differencequotients} is locally bounded.
	Now apply Vitali's theorem to complete the proof of (i).
\end{proof}

In the following we will identify two measurable families of probability measures on the real line, $\{\mu_t,t\in[0,\infty)\}, \{\nu_t,t\in[0,\infty)\}$, if there is a subset $N\subset[0,\infty)$ of Lebesgue measure zero such that $\mu_t=\nu_t$ for all 
$t\in[0,\infty)\backslash N$.

\begin{theorem}\label{T:chordalloewnerequation}
	If $\{f(t;\cdot),t\in[0,\infty)\}$ is any chordal Loewner family, 	then there is a unique measurable family 	$\{\mu_t,t\in[0,\infty)\}$ of probability measures on the real line
	and a subset $N\subset[0,\infty)$ of Lebesgue measure 	zero such that 
	\begin{equation}\label{E:chordalloewnerequation}
		\frac{\partial}{\partial t}f(t;z)=-\int_{\mathbb R}
		\frac{\mu_t(dx)}{z-x}\cdot\frac{\partial}{\partial z}f(t;z)
	\end{equation}
	for all $t\in[0,\infty)\backslash N$, and $z\in\mathbb H$.
\end{theorem}

\begin{proof}
	Let $\{B(a,b;\cdot),0\le a\le b<\infty\}$ be the semigroup 	associated to the chordal Loewner family 	$\{f(t;\cdot),t\in[0,\infty)\}$. Then
	\begin{align}
		\frac{f(b;z)-f(a;z)}{b-a}&=\frac{f(a;B(a,b;z))-f(a;z)}{b-a}
		\notag\\
		&=\frac{f(a;B(a,b;z))-f(a;z)}{B(a,b;z)-z}\cdot
		\frac{B(a,b;z)-z}{b-a}.\notag
	\end{align}
	By \eqref{E:absolutesecond}, $B(a,b;z)\to z$, as 
	$b\searrow a$. So
	\[
		\lim_{b\searrow a}\frac{f(a;B(a,b;z))-f(a;z)}{B(a,b;z)-z}
		=\frac{\partial}{\partial z}f(a;z).
	\]
	Since $(\partial/\partial z)f(a;z)\neq0$, the limit 
	\begin{equation}\label{E:limcauchy}
		\lim_{b\searrow a}\frac{B(a,b;z)-z}{b-a}
	\end{equation}
	also exists. Furthermore, since
	\[
		y\left|\frac{B(a,b;iy)-iy}{b-a}\right|=1+o(1),
		\text{as }y\to\infty,
	\]
	each function $[B(a,b;z)-z]/(b-a)$ is the Cauchy transform of a 	probability measure on $\mathbb R$. By 
	\cite[Theorem 2.5]{maassen:1992}, the limit is also the Cauchy 	transform of a probability measure on $\mathbb R$, say 	$\mu_a$.

	For fixed $z\in\mathbb H$, the function $a\in[0,\infty)\to
	G_a(z)\equiv\int_{\mathbb R}\mu_a(dx)/(z-x)$ is measurable. It 	follows that 
	$a\in[0,\infty)\mapsto G_a\in C(\mathbb H;-\mathbb H)$ is also 	measurable if we endow $C(\mathbb H;-\mathbb H)$, the 	space of continuous functions, with the 	topology of uniform convergence on compact subsets and
	consider the Borel $\sigma$-field. By 
	\cite[Theorem 2.5]{maassen:1992}, the map 
	$a\in[0,\infty)\mapsto \mu_a\in{\bf M}_1(\mathbb R)$ is then	measurable if we endow ${\bf M}_1(\mathbb R)$, the space of  	probability measures on the real line, with the topology of 	weak convergence. Finally, the family $\{\mu_t,t\in[0,\infty)\}$ 	is unique because the measures are determined, $a.e. $ in  	$t$, by \eqref{E:chordalloewnerequation}.
\end{proof}

\begin{theorem}\label{T:loewnersemigroupequation}
	In Theorem \ref{T:chordalloewnerequation}, let 
	$\{B(a,b;\cdot),0\le a\le b<\infty\}$ be the semigroup 	associated to the chordal Loewner family.

	(i) On $[a,\infty)$,
	\[
		\frac{\partial}{\partial s}B(a,s;z)=-\int_{\mathbb R}
		\frac{\mu_s(dx)}{z-x}\cdot\frac{\partial}{\partial z}B(a,s;z),
		\quad B(a,a;z)=z.
	\]
	
	(ii) If $b>0$, then on $[0,b]$
	\[
		\frac{\partial}{\partial t}B(t,b;z)=\int_{\mathbb R}
		\frac{\mu_t(dx)}{B(t,b;z)-x},\quad B(b,b;z)=z.
	\]
\end{theorem}

\begin{proof}
	On $[a,\infty)$,
	\[
		\frac{\partial}{\partial s}f(s;z)=\frac{\partial}{\partial s}
		f(a;B(a,s;z))=f_2(a;B(a,s;z))\frac{\partial}{\partial s}B(a,s;z),
	\]
	and 
	\[
		\int_{\mathbb R}\frac{\mu_s(dx)}{z-x}\cdot
		\frac{\partial}{\partial z}f(s;z)=\int_{\mathbb R}
		\frac{\mu_s(dx)}{z-x}\cdot f_2(a;B(a,s;z))
		\frac{\partial}{\partial z}B(a,s;z).
	\]
	By Theorem \ref{T:absolute} the application of the chain rule is 	valid, \cite[Theorem 8.3C]{rosenblum.rovnyak:1994}. Now 
	\eqref{E:chordalloewnerequation} implies (i).

	If $0\le a\le t\le b$, then $B(a,b;z)=B(a,t;B(t,b;z))$. Therefore
	\begin{align}
		0&=\frac{\partial}{\partial t}B(a,b;z)\notag\\
		&=B_2(a,t;B(t,b;z))+B_3(a,t;B(t,b;z))\frac{\partial}{\partial t}
		B(t,b;z).\notag
	\end{align}
	Since 
	\[
		B_2(a,t;w)=-\int_{\mathbb R}\frac{\mu_t(dx)}{w-x}\cdot
		B_3(a,t;w)
	\]
	by (i), we get
	\[
		B_3(a,t;B(t,b;z))\frac{\partial}{\partial t}B(t,b;z)=
		\int_{\mathbb R}\frac{\mu_t(dx)}{B(t,b;z)-x}\cdot
		B_3(a,t;B(t,b;z)).
	\]
	Again by Theorem \ref{T:absolute} the application of the chain 	rule is justified and, together with $B_3(a,t;w)\neq0$, this 	proves (ii). 
\end{proof} 
	
We now want to show that every measurable family $\{\mu_t,t\in[0,\infty)\}$ of probability measures on the real line determines a unique chordal Loewner family $\{f(t;\cdot),t\in[0,\infty)\}$ such that
\[
	\frac{\partial}{\partial t}f(t;z)=-\int_{\mathbb R}
	\frac{\mu_t(dx)}{z-x}\cdot\frac{\partial}{\partial z}f(t;z).
\]


\begin{theorem}\label{T:generateloewnersemigroup}
	Let $\{\mu_t,t\in[0,\infty)\}$ be a measurable family of 	probability measures on the real line. 	There exists a unique 	family of functions $\{B(a,b;\cdot), 0\le a\le 	b<\infty\}$ with these properties:

	(i) For fixed $a,b$, $B(a,b;\cdot)$ is in $\mathfrak R$, 	$\lim_{y\to\infty}iy[iy-B(a,b;iy)]=b-a$, and 
	\begin{equation}\label{E:flow}
		B(a,c;z)=B(a,b;B(b,c;z))
	\end{equation}
	whenever $0\le a\le b\le c$.

	(ii) For fixed $b>0$ and $z\in\mathbb H$, $a\in[0,b]\mapsto 	B(a,b;z)\in\mathbb H$ is absolutely continuous such that 
	\begin{equation}\label{E:initialvalue}
		\frac{\partial}{\partial a}B(a,b;z)=\int_{\mathbb R}
		\frac{\mu_{a}(dx)}{B(a,b;z)-x},
	\end{equation}
	$a.e.$ on $[0,b]$, and $B(b,b;z)=z$.

	Furthermore, if for each $b\in[0,\infty)$ there exists an 	$M\in\mathbb Z^+$ such that $\text{supp }\mu_a\subset	[-M,M]$ for every $a\in[0,b]$, then
	\[
		B(a,b;z)=z-\frac{b-a}{z}+O\left(\frac{1}{|z|^2}\right),\quad 		z\to\infty.
	\]
\end{theorem}

\begin{proof}
	If $B(\cdot,b;z)$ solves the initial value problem 	\eqref{E:initialvalue}, then it also solves the integral equation
	\begin{equation}\label{E:integral}
		B(a,b;z)=z-\int_a^b\left(
		\int_{\mathbb R}\frac{\mu_{s}(dx)}{B(s,b;z)-x}\right)\ 		ds,\quad 0\le a\le b.
	\end{equation}
	Furthermore, continuous solutions of \eqref{E:integral} satisfy 
	\begin{equation}\label{E:imincrease}
		\Im(B(a,b;z))\ge\Im(z),\quad 0\le a\le b,
	\end{equation}
	since $\Im(-\int\mu(dx)/(z-x))>0$ if $\Im(z)>0$. To solve 	\eqref{E:integral} we construct functions $B_0(a),B_1(a),\dots$, 	$a\in[0,b]$ such that $B_0(a)\equiv z$ and 
	\[
		B_{n+1}(a)=z-\int_a^b\left(\int_{\mathbb R}
		\frac{\mu_{s}(dx)}{B_n(s)-x}\right)\ ds,\quad a\in[0,b],
	\]
	$n=0,1,2,\dots$. Then
	\[
		B_1(a)=z-\int_a^b\left(\int_{\mathbb R}
		\frac{\mu_{s}(dx)}{z-x}\right)\ ds,\quad a\in[0,b],
	\]
	satisfies $\Im(B_1(a))\ge\Im(z), a\in[0,b]$ and by induction it 	follows immediately that $\Im(B_n(a))\ge\Im(z)$ for all 	$n\in\mathbb Z^+$ and $a\in[0,b]$.  Furthermore,
	\[
		|B_1(a)-B_0(a)|=\left|\int_a^b\left(\int_{\mathbb R}
		\frac{\mu_{s}(dx)}{z-x}\right)\ ds\right|\le\frac{b-a}{\Im(z)}.
	\]
	Similarly,
	\begin{align}
		|B_2(a)-B_1(a)|&=\left|\int_a^b(B_2(s)-B_1(s))\left(
		\int_{\mathbb R}
		\frac{\mu_{s}(dx)}{(B_2(s)-x)(B_1(s)-x)}\right)\ 		ds\right|\notag\\
		&\le\frac{1}{\Im(z)^2}\int_a^b|B_1(s)-B_0(s)|\ 		ds\le\frac{1}{\Im(z)^3}\cdot\frac{(b-a)^2}{2!}.\notag
	\end{align}
	Continuing inductively we obtain continuous functions 	$B_0(a),B_1(a),B_2(a),\dots$, $a\in[0,b]$, satisfying
	\[
		|B_{n+1}(a)-B_n(a)|
		\le\frac{1}{\Im(z)^{2n+1}}\cdot\frac{(b-a)^{n+1}}{(n+1)!},
		\quad a\in[0,b],
	\]
 	for each $n\in\mathbb N$. The estimates imply that the limit
	\[
		B(a)=\lim_{n\to\infty}B_n(a)
		=z+\sum_{n=0}^{\infty}(B_{n+1}(a)-B_n(a))
	\]
	exists uniformly on $[0,b]$, that $\Im(B(a))\ge\Im(z)$, 	$a\in[0,b]$,  and that $a\in[0,b]\mapsto B(a)\in\mathbb H$ 	satisfies \eqref{E:integral}. To show that $B$ is the unique 	solution of \eqref{E:integral} suppose that $\tilde{B}$ is 	another continuous function on $[0,b]$ satisfying 	\eqref{E:integral}. Then $\Im(\tilde{B}(a))\ge\Im(z)$ and 
	\begin{align}
		|B(a)-\tilde{B}(a)|&=\left|\int_a^b(\tilde{B}(s)-B(s))\left(
		\int_{\mathbb R}\frac{\mu_{s}(dx)}{(B(s)-x)(\tilde{B}(s)-x)}
		\right)\ ds\right|\notag\\
		&\le\frac{1}{\Im(z)^2}\int_a^b|B(s)-\tilde{B}(s)|\ ds.\notag
	\end{align}
	Now Gronwall's inequality implies $B(a)=\tilde{B}(a)$, 	$a\in[0,b]$. To show that $B(a)=B(a,b;z)$ is analytic as a 	function of $z\in\mathbb H$ note first that $B_0(a)\equiv z$ is 	analytic on $\mathbb H$. Suppose now that 	$B_n(a)=B_n(a,b;z)$ is analytic on $\mathbb H$. Since 	$\Im(B_n(a,b;z))\ge\Im(z)$ we get by bounded convergence 	that $z\in\mathbb H\mapsto B_{n+1}(a,b;z)\in\mathbb H$ is 	continuous and, by Fubini's theorem, that for any closed 	triangle $\Delta$ in $\mathbb H$
	\[
		\int_{\Delta}B_{n+1}(a,b;z)\ dz=\int_a^b\left[
		\int_{\mathbb R}\left(\int_{\Delta}\frac{1}{B_n(s,b;z)-x}\ dz
		\right)\ \mu(dx)\right]\ ds=0.
	\]
	It follows by Morera's theorem that $B_{n+1}(a,b;z)$ is analytic 	as a function of $z$ in $\mathbb H$. Since 
	$B_n(a,b;z)\to B(a,b;z)$, $n\to\infty$, uniformly on compact 	subsets of $\mathbb H$, it follows that $B(a,b;z)$ is analytic 	as a function of $z$ in $\mathbb H$.

	We show that $B(a,b;z)$ is univalent on $\mathbb H$ for 	$a\in[0,b]$. This is clear for $a=b$ and to prove it for 	$a\in[0,b)$ suppose that $B(a_0,b;z_1)=B(a_0,b;z_2)$ for some 	$a_0\in[0,b)$ and $z_1,z_2\in\mathbb H$. Note that 
	\begin{align}
		\frac{\partial}{\partial a}&(B(a,b;z_1)-B(a,b;z_2))\notag\\
		&=(B(a,b;z_1)-B(a,b;z_2))\int_{\mathbb R}
		\frac{\mu_{a}(dx)}{(B(a,b;z_1)-x)(B(a,b;z_2)-x)},\notag
	\end{align}
	for almost every $a\in[0,b]$. Hence, setting 	$w(a)=B(a,b;z_1)-B(a,b;z_2)$, we have $w(a_0)=0$ and 
	\[
		\left|\frac{\partial}{\partial a}w(a)\right|
		\le\frac{1}{\Im(z_1)\Im(z_2)}|w(a)|\quad\text{for  a.e. } 		a\in[0,b].
	\] 
	Choose $M>0$ such that $|w(a)|\le M$ for $a\in[a_0,b]$. Then
	\[
		|w(a)|=\left|\int_{a_0}^a\left(\frac{\partial}{\partial s}w(s) 		\right)\ ds\right|\le\frac{M}{\Im(z_1)\Im(z_2)}(a-a_0)
	\]
	and
	\begin{align}
		|w(a)|&\le\int_{a_0}^a\frac{1}{\Im(z_1)\Im(z_2)}|w(s)|\ ds
		\notag\\					
		&=\frac{1}{\Im(z_1)\Im(z_2)}\int_{a_0}^a\left|\int_{a_0}^s
		\left(\frac{\partial}{\partial t}w(t)\right)\ dt\right|\ ds\notag\\
		&\le\frac{M}{(\Im(z_1)\Im(z_2))^2}\cdot\frac{(a-a_0)^2}{2!}.
		\notag
	\end{align}
	Upon iteration we obtain a sequence of estimates which imply 	that $w(a)=0$ for $a\in[a_0,b]$. In particular, $w(0)=0$ implies 	$z_1=z_2$ and this completes the proof that $B(a,b;z)$ is 	univalent on $\mathbb H$. 

	Next, for fixed $0\le a\le c$ and $z\in\mathbb H$ define 	$t\in[0,c]\mapsto u(t)\in\mathbb H$ by 
	\begin{equation}
		u(t)=
		\begin{cases}
			B(t,a;B(a,c;z)), &\text{if $t\in[0,a]$},\\
			B(t,c;z), &\text{if $t\in(a,c].$}
		\end{cases}
	 \end{equation}
	We will show that $u$ satisfies \eqref{E:integral}. Uniqueness 	of the solution to \eqref{E:integral} then implies the flow 	identity \eqref{E:flow}. To show that $u$ solves 	\eqref{E:integral} note first that $u$ is continuous and that it 	solves \eqref{E:integral} for $t\in(a,c]$ by definition. For 	$t\in[0,a]$ we have
	\begin{align}
		u(t)&=B(a,c;z)-\int_t^a\left(\int_{\mathbb R}
		\frac{\mu_{s}(dx)}{B(s,a;B(a,c;z))-x}\right)\ ds\notag\\
		&=z-\int_a^c\left(\int_{\mathbb R}
		\frac{\mu_{s}(dx)}{B(s,c;z)-x}\right)\ ds\notag\\
		&\qquad-\int_t^a\left(\int_{\mathbb R}
		\frac{\mu_{s}(dx)}{B(s,a;B(a,c;z))-x}\right)\ ds\notag\\
		&=z-\int_t^c\left(\int_{\mathbb R}
		\frac{\mu_{s}(dx)}{u(s)-x}\right)\ ds.\notag
	\end{align}
	Hence $u$ solves \eqref{E:integral} on $[0,c]$ and 	\eqref{E:flow} follows.

	Next, by bounded convergence,
	\begin{align}
		iy& \left[iy-B(a,b;iy)\right]\notag\\
		&=\int_a^b\left(\int_{\mathbb R}
		\frac{iy}{B(s,b;iy)-x}\ \mu_{s}(dx)\right)\ ds\to b-a,\notag
	\end{align}
	as $y\to\infty$.

	Finally, if $b\in[0,\infty)$ and $M\in\mathbb Z^+$ such that 	$\text{supp }\mu_a\subset[-M,M]$ for every $a\in[0,b]$, then 	for all $z\in\mathbb H$ such that $\Im(z)>M$
	\begin{equation}\label{E:expand}
		\int_{\mathbb R}\frac{\mu_{a}(dx)}{B(a,b;z)-x}
		=\frac{1}{B(a,b;z)}+\sum_{k=1}^{\infty}
		\frac{m_{a}(k)}{(B(a,b;z))^{k+1}},
	\end{equation}
	where $m_{a}(k)=\int_{\mathbb R}x^k\ \mu_a(dx)$. Since 	$|(\partial/\partial a)B(a,b;z)|\le1/\Im(z)$ for a.e. $a\in[0,b]$ it is 	clear that
	\[
		\lim_{z\to\infty}(B(a,b;z)-z)=0,
	\]
	and so \eqref{E:integral} and \eqref{E:expand} imply
	\[
		B(a,b;z)=z-\int_a^b\frac{ds}{z}+O\left(\frac{1}{|z|^2}\right),
		\quad z\to\infty,
	\]
	which proves the last statement of the theorem.
\end{proof}

By Theorem \ref{T:chordalloewnerequation}, for every chordal Loewner family $\{f(t;\cdot),t\in[0,\infty)\}$ there is a measurable family $\{\mu_t,t\in[0,\infty)\}$ of probability measures on $\mathbb R$ such that 
\begin{equation}\label{E:cle}
	\frac{\partial}{\partial t}f(t;z)=-\int_{\mathbb R}
		\frac{\mu_t(dx)}{z-x}\cdot\frac{\partial}{\partial z}f(t;z).
\end{equation}
Every measurable family of probability measures on $\mathbb R$ arises in this way.

\begin{theorem}\label{T:measurestoclf}
	If $\{\mu_t,t\in[0,\infty)\}$ is any measurable family of 	probability measures on the real line, then there exists a 	unique chordal Loewner family $\{f(t;\cdot),t\in[0,\infty)\}$ such 	that \eqref{E:cle} holds.
\end{theorem}

\begin{proof}
	Construct $\{B(a,b;\cdot),0\le a\le b<\infty\}$ as in Theorem 	\ref{T:generateloewnersemigroup} and set
	\[
		f(t;z)=B(0,t;z),\quad z\in\mathbb H,t\in[0,\infty).
	\]
	By \eqref{E:flow}, $\{f(t;\cdot),t\in[0,\infty)\}$ is a chain in 	$\mathfrak R$ which, by the statement preceding 	\eqref{E:flow} and Theorem \ref{T:representation}, is maximal, 	i.e. $\{f(t;\cdot),t\in[0,\infty)\}$ is a chordal Loewner family. By 	Theorem \ref{T:chordalloewnerequation}
	\[
		\frac{\partial}{\partial t}f(t;z)=-\int_{\mathbb R}
		\frac{\nu_t(dx)}{z-x}\cdot\frac{\partial}{\partial z}f(t;z)
	\]
	for some measurable family $\{\nu_t,t\in[0,\infty)\}$ of 	probability measures on the real line. Since $f(a;z)=f(t;B(t,a;z))$ 	for $t\in[0,a]$,
	\begin{align}
		0&=f_1(t;B(t,a;z))+f_2(t;B(t,a;z))\frac{\partial}{\partial t}
		B(t,a;z)\notag\\
		&=f_1(t;B(t,a;z))+f_2(t;B(t,a;z))\int_{\mathbb R}
		\frac{\mu_t(dx)}{B(t,a;z)-x}.\notag
	\end{align}
	 It follows that $\mu_t=\nu_t$ $a.e.$ in $t$ and therefore 	$f(t;z)$ satisfies \eqref{E:cle} for the given family 	$\{\mu_t,t\in[0,\infty)\}$.

	Suppose that $\{g(t;\cdot),t\in[0,\infty)\}$ is another chordal 	Loewner family that satisfies \eqref{E:cle}. Then, by the 	generalized chain rule and Theorem 	\ref{T:loewnersemigroupequation} (ii),
	\begin{align}
		\frac{\partial}{\partial t}[g(t;B(t,a;z))]&=g_1(t;B(t,a;z))
		+g_2(t;B(t,a;z))\frac{\partial}{\partial t} B(t,a;z)\notag\\
		&=-\int_{\mathbb R}\frac{\mu_t(dx)}{B(t,a;z)-x}\cdot 
		g_2(t;B(t,a;z))\notag\\
		&\qquad+g_2(t;B(t,a;z))\int_{\mathbb R} 		\frac{\mu_t(dx)}{B(t,a;z)-x}\notag\\
		&=0\notag
	\end{align}
	$a.e.$ on $[0,a]$. Therefore 
	\[
		g(t;B(t,a;z))=g(a;B(a,a;z))=g(a;z)
	\]
	for all $t\in[0,a]$. In particular,
	\[
		g(a;z)=g(0;B(0,a;z))=B(0,a;z)=f(a;z).
	\]
\end{proof}


\end{document}